\documentclass[12pt,oneside]{amsart}

\usepackage{amssymb}
\usepackage{amscd}

\title[Birational embeddings with two Galois points]{Birational embeddings of the Hermitian, Suzuki and Ree curves with two Galois points} 
\author{Satoru Fukasawa}

\address{Department of Mathematical Sciences, Faculty of Science, Yamagata University, Kojirakawa-machi 1-4-12, Yamagata 990-8560, Japan} 
\email{s.fukasawa@sci.kj.yamagata-u.ac.jp} 

\thanks{The author was partially supported by JSPS KAKENHI Grant Number 16K05088.}

\newtheorem{theorem}{Theorem}

\begin{document}
\begin{abstract}  
We show that there exists a plane curve of degree $q^3+1$ with two inner Galois points whose smooth model is the Hermitian curve of degree $q+1$, where $q$ is a power of the characteristic $p>0$.  
Similar results hold for the Suzuki and Ree curves respectively. 
\end{abstract}

\maketitle 

\section{Introduction} 
Let $C \subset \Bbb P^2$ be an irreducible plane curve over an algebraically closed field $k$ of characteristic $p \ge 0$ and let $k(C)$ be its function field.  
For a point $P \in \Bbb P^2$, if the function field extension $k(C)/\pi_P^*k(\Bbb P^1)$ induced by the projection $\pi_P$ is Galois, then $P$ is called a Galois point for $C$. 
This notion was introduced by Yoshihara (\cite{miura-yoshihara, yoshihara1}).
When a Galois point $P$ is a smooth point of $C$, $P$ is called an inner Galois point. 
There are not so many examples of plane curves with two inner Galois points (see the Table in \cite{yoshihara-fukasawa}). 
In this note, we give new examples, which update the Table in \cite{yoshihara-fukasawa}. 

Let $p>0$ and let $q \ge 3$ be a power of $p$. 
We consider the curve $H$ defined by
$$X^qZ+XZ^q-Y^{q+1}=0, $$
which is called the Hermitian curve.  
For the natural embedding of $H$ in $\Bbb P^2$ of degree $q+1$, Homma determined the distribution of Galois points (\cite{homma}). 
To produce examples of plane curves with two Galois points, we would like to consider other birational embeddings. 
We show the following. 

\begin{theorem} \label{hermitian} 
For the Hermitian curve $H$ of degree $q+1$, there exists a morphism $\varphi: H \rightarrow \Bbb P^2$ such that 
\begin{itemize}
\item[(a)] the morphism $\varphi: H \rightarrow \varphi (H)$ is birational, 
\item[(b)] the degree of $\varphi(H)$ is $q^3+1$, and  
\item[(c)] there exist exactly two inner Galois points on $\varphi(H)$. 
\end{itemize} 
\end{theorem} 

For the proof, it is important that two subgroups $G_1$ and $G_2$ of the full automorphism group ${\rm Aut}(H)$ of order $q^3$ act on the set $H(\Bbb F_{q^2})$ of all $\Bbb F_{q^2}$-rational points, which  consists of $q^3+1$ points, on the Hermitian curve $H$.  
The automorphism groups of the Suzuki and Ree curves have the similar property. 
We will obtain the following. 

\begin{theorem} \label{suzuki} 
Let $p=2$, $q_0 \ge 2$ a power of $2$, and let $q=2q_0^2$. 
For the Suzuki curve $\hat{C}$, which is the smooth projective model of the affine curve defined by 
$$ x^q+x=y^{2q_0}(y^q+y), $$
there exists a morphism $\varphi: \hat{C} \rightarrow \Bbb P^2$ such that 
\begin{itemize}
\item[(a)] the morphism $\varphi: \hat{C} \rightarrow \varphi (\hat{C})$ is birational, 
\item[(b)] the degree of $\varphi(\hat{C})$ is $q^2+1$, and  
\item[(c)] there exist exactly two inner Galois points on $\varphi(\hat{C})$. 
\end{itemize} 
\end{theorem}

\begin{theorem} \label{ree} 
Let $p=3$, $q_0 \ge 3$ a power of $3$, and let $q=3q_0^2$. 
For the Ree curve $\hat{C}$, which is the smooth projective model of the affine curve defined by 
$$ y_1^q-y_1=x^{q_0}(x^q-x) \ \mbox{ and } \ y_2^q-y_2=x^{q_0}(y_1^q-y_1), $$
there exists a morphism $\varphi: \hat{C} \rightarrow \Bbb P^2$ such that 
\begin{itemize}
\item[(a)] the morphism $\varphi: \hat{C} \rightarrow \varphi (\hat{C})$ is birational, 
\item[(b)] the degree of $\varphi(\hat{C})$ is $q^3+1$, and  
\item[(c)] there exist exactly two inner Galois points on $\varphi(\hat{C})$. 
\end{itemize} 
\end{theorem}

\section{Hermitian curves}
Let $P_1=(1:0:0)$ and $P_2=(0:0:1) \in H$.
We consider the subgroup  
$$G_{1}:=\left\{
\left(\begin{array}{ccc}
1 & a^q & b \\
0 & 1 & a \\
0 & 0 & 1 
\end{array}\right) \ ; \
a \in \Bbb F_{q^2}, \ b^q+b=a^{q+1}
\right\} $$
of ${\rm PGL}(3, k)$, which is of order $q^3$. 
For any $\sigma \in G_{1}$, it follows that $\sigma(H)=H$, $\sigma(P_1)=P_1$ and $\sigma(H(\Bbb F_{q^2})\setminus \{P_1\})=H(\Bbb F_{q^2}) \setminus \{P_1\}$ (see also \cite[pp.~643--644]{hkt}). 
We take $x=X/Z$ and $y=Y/Z$. 
Note that $k(H)=k(x,y)$ and $y^{q^2}-y \in k(H)^{G_1}$. 
Since $[k(x,y):k(y)]=q$ and $[k(y):k(y^{q^2}-y)]=q^2$, it follows that $k(y^{q^2}-y)=k(H)^{G_1}$ and $k(H)^{G_1} \cong k(\Bbb P^1)$. 
Similarly, we define 
$$G_{2}:=\left\{
\left(\begin{array}{ccc}
1 & 0 & 0 \\
c & 1 & 0 \\
d & c^q & 1 
\end{array}\right) \ ; \
c \in \Bbb F_{q^2}, \ d^q+d=c^{q+1}
\right\}. $$
Note that $k(H)^{G_2}=k((y/x)^{q^2}-(y/x))$. 
Then, $H/G_{i}\cong \Bbb P^1$ for $i=1, 2$, $G_{1} \cap G_{2}=\{1\}$, and 
$$ \{P_1\} \cup \{\sigma(P_2) \ | \ \sigma \in G_{1}\}=H(\Bbb F_{q^2})=\{P_2\}\cup\{\tau(P_1) \ | \ \tau  \in G_{2}\}. $$
It follows from \cite[Theorem 1]{fukasawa} that we have a morphism $\varphi: H \rightarrow \Bbb P^2$ such that $\varphi$ is birational onto $\varphi(H)$, $\deg \varphi(H)=q^3+1$ and there exist two inner Galois points. 

To determine the number of inner Galois points on $\varphi(H)$, we consider the image $\varphi(H(\Bbb F_{q^2}))$. 
As in the proof of \cite[Theorem 1]{fukasawa}, $\varphi$ is represented by 
$$ \left(\frac{1}{y^{q^2}-y}:\frac{x^{q^2}}{y^{q^2}-x^{q^2-1}y}:1\right). $$
Then, $\varphi(P_1)=(0:1:0)$, $\varphi(P_2)=(1:0:0)$ and $\varphi(H(\Bbb F_{q^2}))=\varphi(H) \cap \{Z=0\}$. 
Let $P=(\alpha:\beta:1) \in H(\Bbb F_{q^2}) \setminus \{P_1, P_2\}$. 
Then, $y-\beta$ is a local parameter at $P$. 
Let $u=y-\beta$ and $v=(y/x)-(\beta/\alpha)$. 
Note that
$$y^{q^2}-y=u^{q^2}-u, \ (y/x)^{q^2}-(y/x)=v^{q^2}-v, \mbox{ and }
\frac{x^{q^2}(y^{q^2}-y)}{y^{q^2}-x^{q^2-1}y}=\frac{u}{v} \times \frac{u^{q^2-1}-1}{v^{q^2-1}-1}. $$ 
On the other hand, 
$$ \frac{dv}{dy}=\frac{x-y\frac{dx}{dy}}{x^2}=-x^{q-2}. $$
It follows that $\varphi(P)=(-\alpha^{q-2}:1:0)$. 
When $\alpha^q+\alpha \ne 0$, the fiber $\varphi^{-1}(\varphi(P))$ contains at least $q+1$ points (that is, $\varphi(P)$ is a singular point of $\varphi(H)$). 

We consider the case where $\alpha^{q-1}+1=0$. 
Then, $\varphi(P)=(1:\alpha:0)$ and the projection $\pi_{\varphi(P)}$ is represented by 
$$ \left(-\alpha \frac{1}{y^{q^2}-y}+\frac{x^{q^2}}{y^{q^2}-x^{q^2-1}y}:1\right)=\left(\frac{x-\alpha}{y} \times \frac{(x-\alpha)^{q^2-1}y^{q^2-1}-x^{q^2-1}}{(y^{q^2}-1)(y^{q^2-1}-x^{q^2-1})}:1\right). $$
It follows that the ramification index at $P$ is equal to $q$. 
This implies that the intersection multiplicity of $\varphi(H)$ and the tangent line at $\varphi(P)$ is $q+1$. 

Assume that $\varphi(R)$ is inner Galois. 
Then, the associated Galois group $G_{\varphi(R)}$ is of order $q^3$, which is a Sylow $p$-subgroup of ${\rm Aut}(H) \cong {\rm PGU}(3, q)$ (see \cite[pp.~643--644]{hkt}). 
Therefore, there exists $P \in H(\Bbb F_{q^2})$ such that $\sigma (P)=P$ for any $\sigma \in G_{\varphi(R)}$. 
Then, the order of the pull-back of a linear polynomial given by the tangent line at $P$ is $q^3$ or $q^3+1$, and hence, $\varphi^{-1}(\varphi(P))=\{P\}$. 
It follows that $P=P_1$ or $P_2$, and hence, $R=P_1$ or $P_2$.   
The proof of Theorem \ref{hermitian} is completed. 

\section{Suzuki curves} 
See \cite{hansen-stichtenoth}, \cite{gkt} or \cite[Section 12.2]{hkt} for properties of the Suzuki curves. 
We take $x=X/Z$ and $y=Y/Z$. 
Let $p=2$, $q_0$ a power of $2$, $q=2q_0^2$, and let $C \subset \Bbb P^2$ be (the projective closure of) the curve defined by 
$$ x^q+x=y^{2q_0}(y^q+y).$$ 
The smooth model of $C$ is denoted by $\hat{C}$ with normalization $r: \hat{C} \rightarrow C$. 
Let $P_{\infty}=(1:0:0)$ and $P_2=(0:0:1) \in C$. 
It is known that $P_{\infty}$ is a unique singular point of $C$ and $r^{-1}(P_{\infty})$ consists of a unique point $P_1 \in \hat{C}$. 

We consider the subgroup  
$$G_{1}:=\left\{
\left(\begin{array}{ccc}
1 & a^{2q_0} & b \\
0 & 1 & a \\
0 & 0 & 1 
\end{array}\right) \ ; \
a, b \in \Bbb F_{q}
\right\} $$
of ${\rm PGL}(3, k)$, which is of order $q^2$.  
For any $\sigma \in G_{1}$, it follows that $\sigma(P_{\infty})=P_{\infty}$ and $\sigma(C(\Bbb F_{q})\setminus \{P_\infty\})=C(\Bbb F_{q}) \setminus \{P_\infty\}$.  
In particular, there exists an inclusion $G_1 \hookrightarrow {\rm Aut}(\hat{C})$. 
Note that $k(C)=k(x,y)$ and $y^{q}+y \in k(C)^{G_1}$. 
Since $[k(x,y):k(y)]=q$ and $[k(y):k(y^{q}+y)]=q$, it follows that $k(y^{q}+y)=k(C)^{G_1}$ and $k(C)^{G_1} \cong k(\Bbb P^1)$. 
Let $h:=xy+x^{2q_0}+y^{2q_0+2}$ and let $\psi$ be the rational transformation of $\Bbb A^2$ given by 
$$ (x, y) \mapsto (y/h, x/h). $$
Then, $\psi$ induces an involution of $\hat{C}$ and $\psi(P_1)=P_2$. 
Let $G_2:=\psi G_1 \psi \subset {\rm Aut}(\hat{C})$, which is of order $q^2$. 
Note that $k(C)^{G_2}=k((x/h)^q+(x/h))$. 
Then, $\hat{C}/G_{i}\cong \Bbb P^1$ for $i=1, 2$, $G_{1} \cap G_{2}=\{1\}$, and 
$$ \{P_1\} \cup \{\sigma(P_2) \ | \ \sigma \in G_{1}\}=\hat{C}(\Bbb F_{q})=\{P_2\}\cup\{\tau(P_1) \ | \ \tau  \in G_{2}\}. $$
It follows from \cite[Theorem 1]{fukasawa} that we have a morphism $\varphi: \hat{C} \rightarrow \Bbb P^2$ such that $\varphi$ is birational onto $\varphi(\hat{C})$, $\deg \varphi(\hat{C})=q^2+1$ and there exist two inner Galois points. 

To determine the number of inner Galois points on $\varphi(\hat{C})$, we consider the image $\varphi(\hat{C}(\Bbb F_{q}))$. 
As in the proof of \cite[Theorem 1]{fukasawa}, $\varphi$ is represented by 
$$ \left(\frac{1}{y^{q}+y}:\frac{h^{q}}{x^{q}+h^{q-1}x}:1\right). $$
Then, $\varphi(P_1)=(0:1:0)$, $\varphi(P_2)=(1:0:0)$ and $\varphi(\hat{C}(\Bbb F_{q}))=\varphi(\hat{C}) \cap \{Z=0\}$. 
Let $P=(\alpha:\beta:1) \in C(\Bbb F_{q}) \setminus \{P_\infty, P_2\}$. 
Then, $y+\beta$ is a local parameter at $P$. 
Let $u=y+\beta$ and $v=(x/h)+(\alpha/h(\alpha, \beta))$. 
Note that
$$y^{q}+y=u^{q}+u, \ (x/h)^{q}+(x/h)=v^{q}+v, \mbox{ and }
\frac{h^q(y^{q}+y)}{x^{q}+h^{q-1}x}=\frac{u}{v} \times \frac{u^{q-1}+1}{v^{q-1}+1}. $$ 
On the other hand, 
$$ \frac{dv}{dy}=\frac{y^{2q_0}h+x(y^{2q_0+1}+x)}{h^2}=\frac{h^{2q_0}}{h^2}=h^{2q_0-2}, $$
using the conditions $x^2=h^{2q_0}+x^{2q_0}y^{2q_0}+y^{4q_0+2}$ and $x^{2q_0}=h+xy+y^{2q_0+2}$. 
It follows that $\varphi(P)=(h^{2q_0-2}(\alpha, \beta):1:0)$. 
Note that $h^{q_0-1}(\gamma^{2q_0+1}, 0)=h^{q_0-1}(0, \gamma)=(\gamma^{2q_0+2})^{q_0-1}=1/\gamma$, for any $\gamma \in \Bbb F_q \setminus \{0\}$. 
It follows that the set $\varphi(C(\Bbb F_q)\setminus \{P_{\infty}, P_2\})$ coincides with $\{Z=0\}(\Bbb F_q) \setminus \{(0:1:0), (1:0:0)\}$, and the fiber $\varphi^{-1}(\varphi(P))$ contains at least two points (that is, $\varphi(P)$ is a singular point of $\varphi(\hat{C})$) for any $P \in C(\Bbb F_q)\setminus \{P_{\infty}, P_2\}$. 

Assume that $\varphi(R)$ is inner Galois. 
Then, the associated Galois group $G_{\varphi(R)}$ is of order $q^2$, which is a Sylow $2$-subgroup of the Suzuki group ${\rm Sz}(q)$ (see \cite[p.~564]{hkt}). 
Therefore, there exists $P \in \hat{C}(\Bbb F_{q})$ such that $\sigma (P)=P$ for any $\sigma \in G_{\varphi(R)}$. 
Similar to the proof of Theorem \ref{hermitian}(3), $R=P_1$ or $P_2$.   
The proof of Theorem \ref{suzuki} is completed. 

\section{Ree curves} 
See \cite{pedersen} or \cite[Section 12.4]{hkt} for properties of the Ree curves. 
Let $p=3$, $q_0$ a power of $3$, $q=3q_0^2$, and let $C \subset \Bbb P^3$ be (the projective closure of) the space curve defined by 
$$ y_1^q-y_1=x^{q_0}(x^q-x) \ \mbox{ and } \ y_2^q-y_2=x^{q_0}(y_1^q-y_1),  $$
where $(x, y_1, y_2)$ and $(X:Y_1:Y_2:Z)$ are systems of affine and homogeneous coordinates of $\Bbb A^3$ and $\Bbb P^3$ respectively. 
The smooth model of $C$ is denoted by $\hat{C}$ with normalization $r: \hat{C} \rightarrow C$. 
Let $P_{\infty}=(0:0:1:0)$ and $P_2=(0:0:0:1) \in C$. 
It is known that $P_{\infty}$ is a unique singular point of $C$ and $r^{-1}(P_{\infty})$ consists of a unique point $P_1 \in \hat{C}$. 

We consider the subgroup  
$$G_{1}:=\left\{
\left(\begin{array}{cccc}
1 & 0 & 0 & a \\
a^{q_0} & 1 & 0 & b \\
a^{2q_0} & -a^{q_0} & 1 & c \\
0  & 0 & 0 & 1 
\end{array}\right) \ ; \
a, b, c \in \Bbb F_{q}
\right\} $$
of ${\rm PGL}(4, k)$, which is of order $q^3$.  
For any $\sigma \in G_{1}$, it follows that $\sigma(P_{\infty})=P_{\infty}$ and $\sigma(C(\Bbb F_{q})\setminus \{P_{\infty}\})=C(\Bbb F_{q}) \setminus \{P_\infty\}$.  
In particular, there exists an inclusion $G_1 \hookrightarrow {\rm Aut}(\hat{C})$. 
Note that $k(C)=k(x,y_1, y_2)$ and $x^{q}-x \in k(C)^{G_1}$. 
Since $[k(x,y_1, y_2):k(x)]=q^2$ and $[k(x):k(x^{q}-x)]=q$, it follows that $k(x^{q}-x)=k(C)^{G_1}$ and $k(C)^{G_1} \cong k(\Bbb P^1)$. 
Let $\psi$ be the involution of $\hat{C}$ induced by 
$$ (x, y_1, y_2) \mapsto (w_6/w_8, w_{10}/w_8, w_9/w_8), $$
as in \cite[p.126]{pedersen} or \cite[p. 577]{hkt} (see also \cite{eid, eid-duursma}). 
It follows that $\psi(P_1)=P_2$. 
Let $G_2:=\psi G_1 \psi \subset {\rm Aut}(\hat{C})$, which is of order $q^3$. 
Note that $k(C)^{G_2}=k((w_6/w_8)^q-(w_6/w_8))$. 
Then, $\hat{C}/G_{i}\cong \Bbb P^1$ for $i=1, 2$, $G_{1} \cap G_{2}=\{1\}$, and 
$$ \{P_1\} \cup \{\sigma(P_2) \ | \ \sigma \in G_{1}\}=\hat{C}(\Bbb F_{q})=\{P_2\}\cup\{\tau(P_1) \ | \ \tau  \in G_{2}\}. $$
It follows from \cite[Theorem 1]{fukasawa} that we have a morphism $\varphi: \hat{C} \rightarrow \Bbb P^2$ such that $\varphi$ is birational onto $\varphi(\hat{C})$, $\deg \varphi(\hat{C})=q^3+1$ and there exist two inner Galois points. 

To determine the number of inner Galois points on $\varphi(\hat{C})$, we consider the image $\varphi(\hat{C}(\Bbb F_{q}))$. 
As in the proof of \cite[Theorem 1]{fukasawa}, $\varphi$ is represented by 
$$ \left(\frac{1}{x^{q}-x}:\frac{w_8^{q}}{w_6^{q}-w_8^{q-1}w_6}:1\right). $$
Then, $\varphi(P_1)=(0:1:0)$, $\varphi(P_2)=(1:0:0)$ and $\varphi(\hat{C}(\Bbb F_{q}))=\varphi(\hat{C}) \cap \{Z=0\}$. 
Let $P=(\alpha:\beta:\gamma:1) \in C(\Bbb F_{q}) \setminus \{P_\infty, P_2\}$. 
Then, $x-\alpha$ is a local parameter at $P$. 
Let $u=x-\alpha$ and $A=(w_6/w_8)-(w_6/w_8)(\alpha, \beta)$. 
Note that
$$x^{q}-x=u^{q}-u, \ (w_6/w_8)^{q}-(w_6/w_8)=A^{q}-A, \mbox{ and }
\frac{w_8^q(x^{q}-x)}{w_6^{q}-w_8^{q-1}w_6}=\frac{u}{A} \times \frac{u^{q-1}-1}{A^{q-1}-1}. $$ 
On the other hand, 
$$ \frac{dA}{dx}=\frac{w_4^{3q_0}w_8-w_6w_7^{3q_0}}{w_8^2}=\frac{w_8^{3q_0}}{w_8^2}=w_8^{3q_0-2}. $$
It follows that $\varphi(P)=(w_8^{3q_0-2}(\alpha, \beta, \gamma):1:0)$. 
Note that $w_8^{3q_0-2}(\delta^{-1}, 0, 0)=-\delta^2$, $w_8^{3q_0-2}(0, \delta, 0)=\delta^{3q_0-3}$ and $w_8^{3q_0-2}(0, \delta^{-q_0-1},0)=\delta^2$, for any $\delta \in \Bbb F_q \setminus \{0\}$. 
It follows from $\sqrt{-1} \not\in \Bbb F_q$ that the set $\varphi(C(\Bbb F_q)\setminus \{P_{\infty}, P_2\})$ coincides with $\{Z=0\}(\Bbb F_q) \setminus \{(0:1:0), (1:0:0)\}$, and the fiber $\varphi^{-1}(\varphi(P))$ contains at least two points (that is, $\varphi(P)$ is a singular point of $\varphi(\hat{C})$) for any $P \in C(\Bbb F_q)\setminus \{P_{\infty}, P_2\}$. 

Assume that $\varphi(R)$ is inner Galois. 
Then, the associated Galois group $G_{\varphi(R)}$ is of order $q^3$, which is a Sylow $3$-subgroup of the Ree group ${\rm Ree}(q)$ (see \cite[p.~575]{hkt}). 
Similar to the proof of Theorem \ref{hermitian}(3),  $R=P_1$ or $P_2$.   
The proof of Theorem \ref{ree} is completed.


\begin{thebibliography}{20} 
\bibitem{eid} A. Eid, Invariant embeddings of the Deligne--Lusztig curves with applications, PhD thesis, University of Illinois at Urbana-Champaign, 2013. 
\bibitem{eid-duursma} A. Eid and I. Duursma, Smooth embeddings for the Suzuki and Ree curves, in {\it Algorithmic Arithmetic, Geometry, and Coding Theory}, pp. 251--291, Contemp. Math. {\bf 637}, Amer. Math. Soc., Providence, RI, 2015. 
\bibitem{fukasawa} S. Fukasawa, A birational embedding of an algebraic curve into a projective plane with two Galois points, preprint, arXiv:1611.03953. 
\bibitem{gkt} M. Giulietti, G. Korchm\'{a}ros and F. Torres, Quotient curves of the Suzuki curve, Acta Arith. {\bf 122} (2006), 245--274. 
\bibitem{hansen-stichtenoth} J. P. Hansen and H. Stichtenoth, Group codes on certain algebraic curves with many rational points, Appl. Algebra Eng. Comm. Comput. {\bf 1} (1990), 67--77. 
\bibitem{hkt} J. W. P. Hirschfeld, G. Korchm\'{a}ros and F. Torres, {\it Algebraic curves over a finite field}, Princeton Univ. Press, Princeton, 2008. 
\bibitem{homma} M. Homma, Galois points for a Hermitian curve, Comm. Algebra {\bf 34} (2006), 4503--4511. 
\bibitem{miura-yoshihara} K. Miura and H. Yoshihara, Field theory for function fields of plane quartic curves, J. Algebra {\bf 226} (2000), 283--294. 
\bibitem{pedersen}  J. P. Pedersen, A function field related to the Ree group, in {\it Coding theory and Algebraic Geometry}, pp. 122--131, Lecture Notes in Math. {\bf 1518}, Springer, Berlin, 1992.  
\bibitem{yoshihara1} H. Yoshihara, Function field theory of plane curves by dual curves, J. Algebra {\bf 239} (2001), 340--355. 
\bibitem{yoshihara-fukasawa} H. Yoshihara and S. Fukasawa, List of problems, available at: \\ http://hyoshihara.web.fc2.com/openquestion.html
\end{thebibliography}
\end{document}